\documentclass[journal,twoside,print]{IEEEtran}
\usepackage{amsmath,amssymb,amsfonts}
\usepackage{algorithmic}
\usepackage{graphicx}
\usepackage{textcomp}
\usepackage{booktabs}

\usepackage{enumitem}
\usepackage{caption}
\usepackage{subcaption}
\DeclareMathOperator{\mvDE}    {mvDE}
\DeclareMathOperator{\mvDEG}    {mvDE_G}
\newcommand{\Na}{\mathbb{N}}
\newcommand{\R}{\mathbb{R}} % symbol for real numbers
\newcommand{\Id}{\mathbb{I}} % symbol for real numbers
\newcommand{\A}{\mathbf{A}}
\newcommand{\N}{\mathbb{N}}
\newcommand{\map}[3]{ #1 \colon #2 \rightarrow #3}
\newcommand{\Y}{\mathbf{Y}}

\newcommand{\card}[1]{\lvert#1\rvert}
\newcommand{\V}{\mathcal{V}} 
\newcommand{\E}{\mathcal{E}} 
\newcommand{\G}{{G}}
\newcommand{\X}{\mathbf{X}} 
\usepackage{cite}

\begin{document}
\title{Graph-Based Multivariate Multiscale Dispersion Entropy: Efficient Implementation and Applications to Real-World Network Data}
\author{John Stewart Fabila-Carrasco, Chao Tan, and Javier Escudero
\thanks{J.S.F.C is with
		the School of Informatics, University of Edinburgh, Edinburgh, EH8 9LE, UK
		(e-mail John.Fabila@ed.ac.uk) }
\thanks{C.T. is  with the School of Electrical and
	Information Engineering
	Tianjin University
	Tianjin, 300072, China
	(e-mail: tanchao@tju.edu.cn). }		
\thanks{J.E is with
the School of Engineering, Institute for Imaging, Data and Communications, University of Edinburgh, Edinburgh, EH9 3FB, UK
(e-mail J.E. Javier.Escudero@ed.ac.uk). }
\thanks{The work of J.S.F.C. and J.E. was supported by Leverhulme Trust through a Research Project under Grant
	RPG-2020-158. For the purpose of open access, the author has applied a Creative Commons Attribution (CC BY) licence to any Author Accepted Manuscript version arising from this submission. }
}

\maketitle

\begin{abstract}
	We introduce Multivariate Multiscale Graph-based Dispersion Entropy (\(\mvDEG\)), a novel, computationally efficient method for analysing multivariate time series data in graph and complex network frameworks, and demonstrate its application in real-world data. \(\mvDEG\) effectively combines temporal dynamics with topological relationships, offering enhanced analysis compared to traditional nonlinear entropy methods. Its efficacy is established through testing on synthetic signals, such as uncorrelated and correlated noise, showcasing its adeptness in discerning various levels of dependency and complexity. 	
	The robustness of \(\mvDEG\) is further validated with real-world datasets, effectively differentiating various two-phase flow regimes and capturing distinct dynamics in weather data analysis. 
	An important advancement of \(\mvDEG\) is its computational efficiency. Our optimized algorithm displays a computational time that grows linearly with the number of vertices or nodes, in contrast to the exponential growth observed in classical methods. This efficiency is achieved through refined matrix power calculations that exploit matrix and Kronecker product properties, making our method faster than the state of the art. The significant acceleration in computational time positions \(\mvDEG\) as a transformative tool for extensive and real-time applications, setting a new benchmark in the analysis of time series recorded at distributed locations and opening avenues for innovative applications.

\end{abstract}


\section{Introduction}
\label{sec:introduction}
Real-world data, particularly from systems like industrial and weather phenomena, is often characterized by dynamic and complex behaviours that exhibit transient, chaotic, and nonlinear dynamics over time, accompanied by a mixture of randomness and uniformity in spatial distribution~\cite{jhaveri2021fault, Tan2013}. These systems are fundamental to various phenomena across environmental, technological, and scientific domains, influencing everything from water distribution~\cite{wei2019optimal} to advanced manufacturing processes~\cite{chai2019enhanced}. Gaining a deep understanding of the underlying mechanisms of these systems is crucial. Despite substantial advancements in the field, a notable challenge persists in the development of sophisticated mathematical tools. These tools are essential not only for an accurate analysis of specific applications like two-phase flow and weather dynamics but also for understanding and analysing general network behaviours across a wide array of domains, including any form of networked or graph-based data, addressing both academic research needs and practical applications~\cite{zhou2020edm}.

Multivariate entropy techniques are essential tools in the analysis of data consisting of multiple time series, or channels~\cite{ahmed2011multivariate, wang2021variational}. These techniques offer a perspective of complexity that extends beyond the scope of univariate analysis~\cite{Ahmed2011, Rostaghi2016}. Among these, classical Multivariate Dispersion Entropy (\(\mvDE\)) is particularly noteworthy for its robust approach in quantifying the complexity and dynamics of multivariate data~\cite{Azami2019}. Its wide-ranging applications, from biomedicine to fault diagnosis, demonstrate its significance in various research fields~\cite{azami2018coarse}. In comparison to other entropy methods such as the Multivariate Sample Entropy~\cite{Ahmed2011}, \(\mvDE\) has better performance, more stability, and is particularly effective for shorter time series~\cite{Azami2019}. 

The landscape of data analysis is rapidly evolving, increasingly pivoting towards the utilization of graph-based structures and complex networks~\cite{fan2023graph, ortega2018graph}. This shift is driven by technological advancements that have vastly increased our capacity to collect and analyse data across diverse fields, ranging from environmental monitoring to public health systems~\cite{boccaletti2006complex, prasse2020network}. Graph-based methods are gaining prominence as they offer a powerful framework to uncover complex relationships and dependencies within this data, which traditional tabular data formats and analysis techniques might not fully capture. These methods leverage the unique topological properties of graph-structured data, enabling a more nuanced understanding of the intricate interplay within  these varied systems~\cite{ortega2018graph}. The development of Permutation Entropy and Dispersion Entropy tailored for graph data marks a significant advancement, facilitating analyses that incorporate these crucial topological dimensions~\cite{Fabila23, Fabila22, fan2023graph}. 

Moreover, the approach we have developed for Permutation and Dispersion Entropy on graph structures has influenced subsequent methodologies, as evidenced by recent advancements like the extension of Bubble Entropy to graph signals~\cite{dong2024information}, directly utilising our foundational work. This not only demonstrates the robustness of our methods but also indicates that our approach is pioneering paths for further innovations in the analysis of network data. However, while these methods provide valuable insights into the \textit{topological dimension} of the data, they often do not concurrently consider the \textit{temporal dimension}. This limitation is crucial, especially in applications where understanding the interplay between topological (or spatial) structure and temporal dynamics is essential for a comprehensive analysis, such as climate science~\cite{qiu2017time} or industrial applications~\cite{dang2019novel}.

For example, in two-phase flow analysis, various entropy metrics have been key to understanding the system dynamics. Sample Entropy has been used to measure regularity in fluid dynamics time-series~\cite{azami2018coarse}. Permutation Entropy helps in detecting changes in flow patterns due to its sensitivity to temporal variations~\cite{tan2023combinational}. Despite their utility, these metrics often do not fully capture the combined spatial and temporal complexities of two-phase flow systems. This limitation highlights the need for more advanced tools like \(\mvDEG\), which are capable of integrating both topological and temporal data for a more comprehensive analysis.

\subsection{Contributions}
We introduce the Multiscale Multivariate Graph-based Dispersion Entropy (\(\mvDEG\)), an advanced method for analysing data on complex networks. Differing from classical Dispersion Entropy~\cite{Rostaghi2016}, which focuses on temporal information, and Dispersion Entropy for Graphs~\cite{Fabila23}, which emphasizes graph structure, \(\mvDEG\) integrates both temporal and topological dimensions for a more comprehensive analysis.

Furthermore, the proposed algorithmic implementation notably increases computational efficiency, making \(\mvDEG\) suitable for processing both short and long time series. Our optimized algorithm displays a computational time that grows linearly with the number of vertices or nodes, in contrast to the exponential growth observed in classical methods. \(\mvDEG\) demonstrates improved performance over classical \(\mvDE\) in various test scenarios, including synthetic and real-world datasets. The efficiency of \(\mvDEG\) can be extended to other graph-based entropy metrics, notably, the principles underlying \(\mvDEG\) could be seamlessly adapted to improve the efficiency of the Multivariate Multiscale Permutation Entropy~\cite{fabila2022multivariate}. Specifically, algorithms that use Cartesian Graph products combined with our adjacency matrix approach, as detailed in previous works on Permutation Entropy~\cite{Fabila22} and extended to Bubble Entropy~\cite{dong2024information}, will find our method straightforwardly applicable. This adaptability ensures that our efficient approach can significantly benefit a wide range of graph-based entropy calculations, fostering advancements across various fields.

In applying \(\mvDEG\) to real-world data including two-phase flow systems and weather phenomena, we have successfully delineated distinct entropy profiles for various complex patterns, providing acute insights into their dynamics. The precision of the method is especially notable in distinguishing different flow regimes and weather patterns, revealing clear entropy demarcations at lower scales. These results not only demonstrate the robustness of \(\mvDEG\) in handling intricate data from industrial processes and environmental monitoring but also underscore its potential as an indispensable tool in advanced time series analysis.

\subsection{Graph Theory Notation}

A \emph{graph} $G$ is defined as $G = (\V,\E)$, where $\V=\{1,2,3,\dots, n\}$ is the set of vertices, and $\E \subset \{(i,j): i,j\in\V\}$ is a set of edges~\cite{bondy1982graph}. The adjacency matrix $\A$, a $N \times N$ matrix, represents the connections between vertices, with $\A_{i j}= 1$ if $(i,j)\in \E$ and $0$ otherwise. For weighted graphs, \(\A_{ij}\) represents the weight of the connection between \(i\) and \(j\). Our approach is applicable to any weighted graph, whether directed or undirected, by appropriately adjusting the adjacency matrix \(\A\). 

The \emph{Cartesian product} of two graphs $\G= (\V,\E)$ and $\G'=(\V',\E')$, denoted by $\G\square \G'$, is a graph defined as:
\begin{enumerate}[wide, labelwidth=!, labelindent=0pt]
	\item $\V(\G\square \G')=\V\times \V'=\{(v,v') \mid v\in\V \text{ and } v'\in\V'\}$.
	\item $(v,v')$ and $(u,u')$ are adjacent in $\G\square \G'$ if and only if either: $v=u$ and $v'$ is adjacent to $u'$ in $\G'$, or $v'=u'$ and $v$ is adjacent to $u$ in $\G$. 
\end{enumerate}

%A \emph{graph signal} is a function $\X: \V \rightarrow \R$ that assigns a real value to each vertex in the graph. This signal is represented as an $n$-dimensional column vector.

\section{Multivariate Multiscale Graph-based Dispersion Entropy ($\mvDEG$)} 
The algorithm $\mvDEG$ enhances dispersion entropy by integrating topological and temporal data dimensions. It builds on the topological approach in~\cite{Fabila23} and the temporal focus of~\cite{Azami2019}, offering a more thorough analysis. However, the initial computation of large matrix powers in $\mvDEG$ presents a computational challenge, which is addressed in Sec.~\ref{sec:fast} through an efficient implementation strategy. The algorithm consists of two main steps: 1) a coarse-graining process, and 2) the calculation of $\mvDEG$ at each scale $\tau$.

Consider a multivariate signal $\mathbf{X}=\left\{x_{k, i}\right\}_{k=1,2, \ldots, p}^{i=1, \ldots, N}$, where $p$ is the number of channels and $N$ is the number of observations per channel. 

\subsection{Coarse-Graining Process for Multivariate Signals}
For each channel, the original signal is divided into non-overlapping segments of length $\tau$, named scale factor. Next, for each channel, the average of each segment is calculated to derive the coarse-grained signals as follows: 
\begin{equation*}
z_{k, i}^{(\tau)}=\frac{1}{\tau} \sum_{b=(i-1) \tau+1}^{i \tau} x_{k, b}, 1 \leq i \leq\left\lfloor\frac{N}{\tau}\right\rfloor=L, 1 \leq k \leq p\;;
\end{equation*}
where $L$ denotes the length of the coarse-grained signal. While this study uses a straightforward coarse-graining approach, various alternative coarse-graining methods are also explored in the literature~\cite{azami2018coarse}. In the second step, $\mvDEG$ is calculated for each coarse-grained signal.

\subsection{Graph-Based Multivariate Dispersion Entropy}
To analyse the inter-channel interactions within $\mathbf{X}$, we utilize an adjacency matrix $I_p \in \mathbb{R}^{p\times p}$. This matrix, $I_p$, encapsulates the connectivity or interaction patterns between the channels. It can be predefined based on system-specific knowledge, or fully connected graph for simplicity. Alternatively, $I_p$ can be inferred from the data itself, employing various graph learning techniques such as in~\cite{xia2021graph}. This flexibility in defining $I_p$ allows for tailored analysis of the multivariate signal, taking into account the unique interaction dynamics present in different datasets or systems.
The $\mvDEG$ is defined through the following steps:

\begin{enumerate}[wide, labelwidth=!, labelindent=0pt]
	\item \textbf{Embedding Matrix Construction:} Given an embedding dimension $2 \leq m \in \Na$ and a class number $c \in \Na$, the embedding matrix $\textbf{Y} \in \R^{N \times m}$ is formed as $\textbf{Y} = [\textbf{y}_0, \textbf{y}_1, \cdots, \textbf{y}_{m-1}]$. Each column $\textbf{y}_k = D\A^{k} \textbf{v}$, where $\A$ is the adjacency matrix of the graph ${\overrightarrow{P}_N} \square I_p$, $D$ is a diagonal normalization matrix, and $\textbf{v}$ is a vectorized form of the multivariate signal $\mathbf{X}$.
	
	\item \textbf{Mapping to Classes:} A mapping function $\map{M}{\R}{\N_c}$, where $\N_c = \{1, 2, \dots, c\}$, is applied element-wise to $\Y$ to classify each entry into one of $c$ classes. This results in a matrix $M(\Y) \in \N_c^{N \times m}$.
	
	\item \textbf{Dispersion Patterns and Relative Frequencies:} Each row of $M(\Y)$, termed an embedding vector, corresponds to a unique dispersion pattern. The relative frequency of each pattern $\pi \in \Pi$ is calculated across the dataset.
	
	\item \textbf{Entropy Computation:} The $\mvDEG$ is the normalized Shannon entropy of the dispersion patterns, computed as:
	\begin{equation*}
		\mvDEG(\X, m, L, c) = -\dfrac{1}{\log(c^m)}\sum_{\pi \in \Pi} p(\pi) \ln p(\pi).
	\end{equation*}
\end{enumerate}

This approach effectively captures the complex interactions within multivariate data, leveraging both the temporal sequence and the underlying graph structure. For embedding matrix construction details, refer to~\cite{Fabila23}. The primary computational challenge involves computing powers of the large $Np \times Np$ matrix $\mathbf{A}_{{\overrightarrow{P}_N} \square I_p}$, representing the graph's adjacency matrix. We address this challenge and its solutions in Sec.~\ref{sec:fast}.

\section{Efficient Algorithm Implementation}\label{sec:fast}
Utilizing matrix properties and the Kronecker product~\cite{graham2018kronecker}, this section introduces an efficient $\mvDEG$ implementation to overcome computational challenges of prior graph signal entropy metrics~\cite{Fabila22}.

\begin{table}[h]
	\centering
	\begin{tabular}{@{}lccccc@{}}
		\toprule
		\textbf{Graph} & Model & \textbf{Notation} & $\card{\V}$ & \textbf{Adjacency} \\ 
		\midrule
		Directed Path& Time & ${\overrightarrow{P}_N}$ & $N$ & $\mathbf{A}_{\overrightarrow{P}_N}$ \\
		Interaction Matrix &  Topology & ${I_p}$ & $p$ & $\mathbf{A}_{I_p}$ \\
		\bottomrule
	\end{tabular}
	\caption{Overview of the graphs and matrices.}
\end{table}
The adjacency matrix of ${\overrightarrow{P}_N} \square I_p$ is given by:
\begin{equation}\label{eq:product}
	\mathbf{A}_{{\overrightarrow{P}_N} \square I_p} = \mathbf{A}_{\overrightarrow{P}_N} \otimes \Id_p + \Id_N \otimes \mathbf{A}_{I_p}
\end{equation}
where $\otimes$ denotes the Kronecker product and $\Id_N$ is the $N \times N$ identity matrix. Let $X = \mathbf{A}_{\overrightarrow{P}_N} \otimes \Id_p$ and $Y = \Id_N \otimes \mathbf{A}_{I_p}$. We demonstrate that \(B\) and \(C\) commute:
\begin{align*}
	BC &= (\mathbf{A}_{\overrightarrow{P}_N} \otimes \Id_p)(\Id_N \otimes \mathbf{A}_{I_p})= \mathbf{A}_{\overrightarrow{P}_N} \Id_N \otimes \Id_p \mathbf{A}_{I_p} \\
	& = \mathbf{A}_{\overrightarrow{P}_N} \otimes \mathbf{A}_{I_p} = \Id_N \mathbf{A}_{\overrightarrow{P}_N} \otimes \mathbf{A}_{I_p} \Id_p\\
	&= (\Id_N \otimes \mathbf{A}_{I_p})(\mathbf{A}_{\overrightarrow{P}_N} \otimes \Id_p)=CB\;.
\end{align*}

Consequently, the $m$-th power matrix in Eq.~\ref{eq:product} is:
\begin{align*}
	(\mathbf{A}_{{\overrightarrow{P}_N} \square I_p})^m &= \sum_{k=0}^{m} \binom{m}{k} (\mathbf{A}_{\overrightarrow{P}_N} \otimes \Id_p)^{k} (\Id_N \otimes \mathbf{A}_{I_p})^{m-k}\\
	&= \sum_{k=0}^{m} \binom{m}{k} ((\mathbf{A}_{\overrightarrow{P}_N})^{k} \otimes \Id_p) (\Id_N \otimes (\mathbf{A}_{I_p})^{m-k})\\
	&=\sum_{k=0}^{m} \binom{m}{k} ((\mathbf{A}_{\overrightarrow{P}_N})^{k} \Id_n)\otimes  (\Id_p (\mathbf{A}_{I_p})^{m-k})\\
	&= \sum_{k=0}^{m} \binom{m}{k} (\mathbf{A}_{\overrightarrow{P}_N})^{k} \otimes (\mathbf{A}_{I_p})^{m-k}\;.
\end{align*}

%
%The adjacency matrix of \( \overrightarrow{P}_N \) is defined as:
%\[
%\mathbf{A}_{\overrightarrow{P}_N}(i, j) = 
%\begin{cases} 
%	1 & \text{if } j = i + 1 \text{ and } j \leq N \\
%	0 & \text{otherwise}
%\end{cases}
%\]
%
%The \( k \)-th power of this matrix is:
%\[
%(\mathbf{A}_{\overrightarrow{P}_N})^k(i, j) = 
%\begin{cases} 
%	1 & \text{if } j = i + k \text{ and } j \leq N \\
%	0 & \text{otherwise}
%\end{cases}
%\]

Finally, the adjacency matrix of \( \overrightarrow{P}_N \), \(\mathbf{A}_{\overrightarrow{P}_N}\), and its \( k \)-th power for \( 1 \leq k \leq N \) are computed as:
\[
\mathbf{A}_{\overrightarrow{P}_N}^k(i, j) = 
\begin{cases} 
	1 & \text{if } j = i + k \text{ and } j \leq N, \\
	0 & \text{otherwise},
\end{cases}
\]
where \(\mathbf{A}_{\overrightarrow{P}_N}(i, j) = \mathbf{A}_{\overrightarrow{P}_N}^1(i, j)\) for \( k = 1 \).

This efficient implementation reduces the computation of the $m$-th power of a $Np \times Np$ matrix to a sum of Kronecker products involving smaller $p \times p$ matrices. This approach effectively addresses the high computational cost associated with graph-based entropy metrics, a major challenge noted in previous studies~\cite{Fabila22, Fabila23}.

\section{Synthetic Signals}

Synthetic signals like multivariate \(1/f\) noise and White Gaussian Noise (WGN) are essential for validating multivariate multiscale entropy techniques. Multivariate \(1/f\) noise, common in natural phenomena, represents complex, correlated systems and contrasts with the high irregularity of WGN~\cite{Azami2019, costa2002multiscale, humeau2016multivariate}. This contrast enables a detailed evaluation of entropy methods, as these signals provide benchmarks to assess how entropy measures handle different complexities and predictabilities in signal patterns, vital for applications~\cite{zhou2020edm}.

%
%The employment of synthetic signals, such as multivariate \(1/f\) noise and White Gaussian Noise (WGN), plays a critical role in the validation and benchmarking of multivariate multiscale entropy techniques. These synthetic signals are integral to entropy analysis, as they allow for a controlled assessment of algorithmic performance against well-understood standards. In the context of entropy analysis, multivariate \(1/f\) noise is especially relevant due to its ubiquity in natural phenomena and its representation of complex, correlated systems~\cite{Azami2019}. Its comparison to WGN, which represents the other end of the spectrum with its high irregularity and uncorrelated nature, provides a comprehensive range for methodological testing~\cite{costa2002multiscale, humeau2016multivariate}. The contrast between the higher complexity but lower irregularity of \(1/f\) noise and the lower complexity but higher irregularity of WGN allows for a nuanced evaluation of entropy methods. These synthetic signals serve as benchmarks to understand how different entropy measures respond to varying levels of complexity and predictability in signal patterns, an aspect crucial for practical applications in fields ranging from biomedical signal processing to industrial data~\cite{zhou2020edm}.

\subsection{Uncorrelated Noise Analysis}
We generate a trivariate time series \(\mathbf{F}(p)  \), where \( p \in \{0,1,2,3\} \) represents the number of channels with \(1/f\) noise, and \( 3-p \) channels with WGN, i.e., \(\mathbf{F}(p) = \{ \underbrace{\text{WGN}, \ldots, \text{WGN}}_{p \text{ times}}, \underbrace{1/f, \ldots, 1/f}_{3-p \text{ times}} \}\;.\)

Here, \( p \) allows a systematic transition from all channels representing independent \(1/f\) noise to all channels being independent WGN. For algorithms on generating these multivariate signals, see~\cite{Azami2019,humeau2016multivariate, costa2002multiscale}.

\begin{figure}
	\centering
	\includegraphics[width=.28\textwidth]{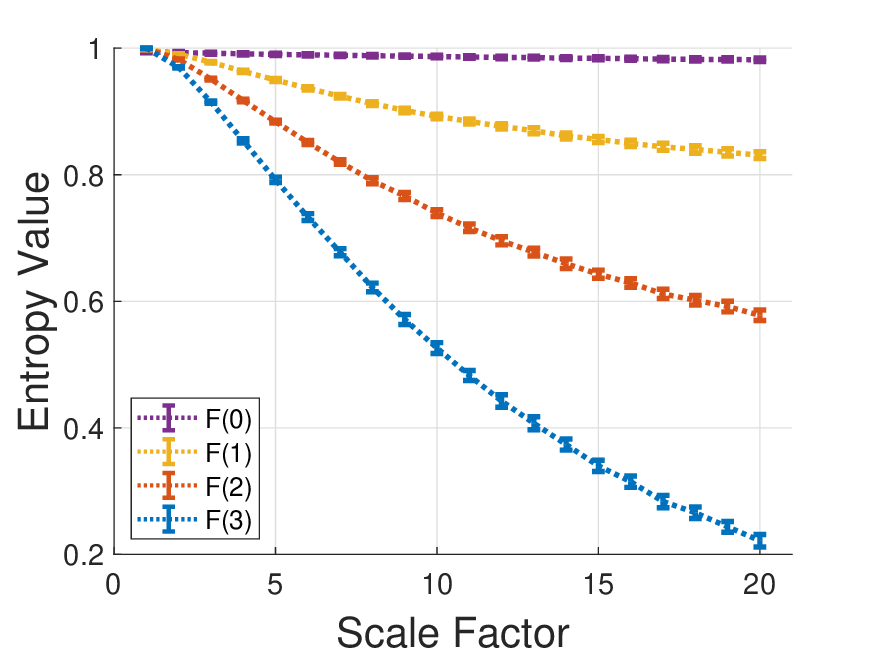}
	\captionsetup{justification=centering}
	\caption{Mean value and standard deviation (SD) for the classical multivariate Dispersion Entropy $\mvDE$.}
	\label{fig:entropy_vs_scale}
\end{figure}

\begin{figure}
	\centering
	\begin{subfigure}{0.24\textwidth}
		\centering
		\includegraphics[width=\textwidth]{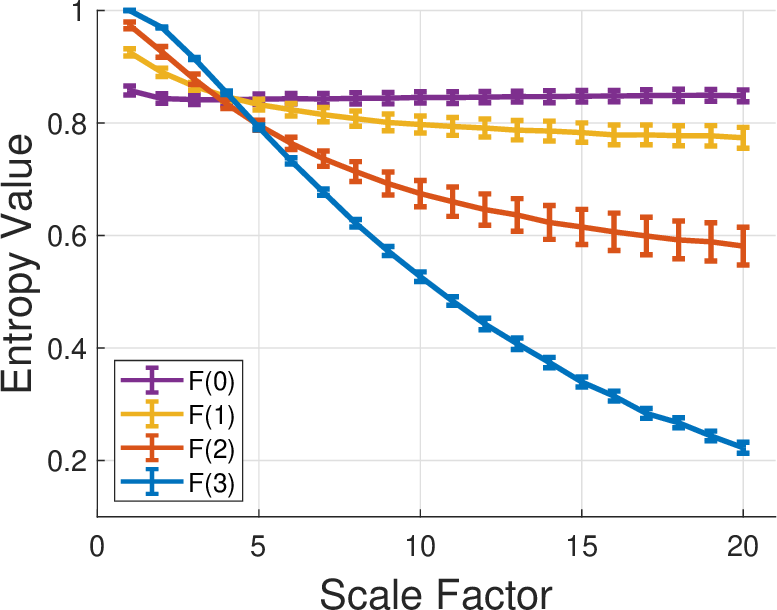}
		\caption{}
		\label{fig:mvde_empty}
	\end{subfigure}
	\begin{subfigure}{0.24\textwidth}
		\centering
		\includegraphics[width=\textwidth]{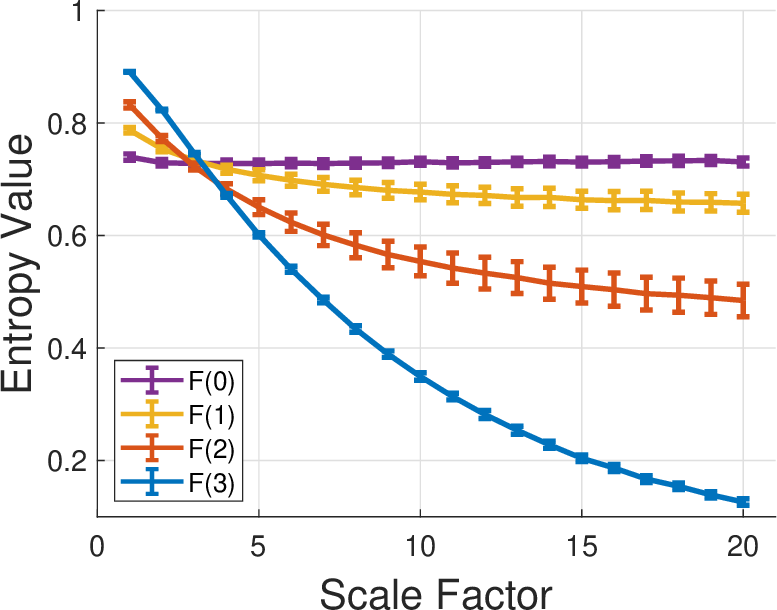}
		\caption{}
		\label{fig:mvdeg_corr}
	\end{subfigure}
	\caption{Entropy value comparison across scale factors for \(\mvDEG\), the mean and SD with (a) identity adjacency graph, indicating no prior inter-channel correlation, and (b) correlation matrix graph, reflecting inherent inter-channel dependencies.}
	\label{fig:mvde_comparison}
\end{figure}

We applied the $\mvDE$ and $\mvDEG$ to $40$ independent realizations of uncorrelated trivariate WGN and \(1/f\) noise, each with $15,000$ sample points. The mean and standard deviations (SD) of the results for $\mvDE$ and $\mvDEG$ are shown for $m=4$ and $c=6$ in both cases in Figs.~\ref{fig:entropy_vs_scale} and \ref{fig:mvde_comparison}, respectively. These \(m\) and \(c\) values are consistent with literature recommendations~\cite{azami2018coarse, Rostaghi2016}.

In both Figs.~\ref{fig:entropy_vs_scale} and \ref{fig:mvde_comparison}, the x-axis represents the Scale Factor, while the y-axis delineates the normalised Entropy Value. The results from both methods are strikingly similar, with the entropy values of trivariate WGN signals being higher than those of other trivariate time series at low scale factors. However, the entropy values for coarse-grained trivariate \(1/f\) noise signals are almost constant along the temporal scale factor, while those for coarse-grained WGN signal monotonically decrease with increasing scale factors. This indicates that multivariate WGN time series contain information primarily at small temporal scale factors, whereas trivariate \(1/f\) noise signals exhibit higher complexity across multiple time series. These observations align with literature findings~\cite{costa2002multiscale}.

In Fig.~\ref{fig:mvde_empty}, we consider the case of a zero adjacency graph ($I_p$ is a $3\times 3$ matrix with all zero entries) since the noise is uncorrelated and theoretically, there is no interaction between time series. The results are similar to $\mvDE$, reflecting the primarily temporal nature of the information. However, the simulated time series exhibit slight correlations in practice. When these correlations are calculated and used as $I_p$ in Fig.~\ref{fig:mvdeg_corr}, the entropy values of $\mvDEG$ are slightly lower, indicating a decrease in complexity due to the estimated correlations. These results are consistent with the notion that \(1/f\) noise is structurally more complex than multivariate WGN.

In single-channel scenarios, multivariate entropy metrics reduce to their univariate counterparts. Thus, classical Dispersion Entropy~\cite{Rostaghi2016} and Dispersion Entropy for Graph Signals~\cite{Fabila23} yield identical results. In multivariate cases, while both \(\mvDE\) and \(\mvDEG\) produce different entropy values, they similarly characterize dynamics. However, \(\mvDEG\) distinguishes itself with its computational efficiency. Unlike \(\mvDE\), where the number of patterns increases exponentially with parameters \(m\) and \(p\) without considering the graph structure, \(\mvDEG\) incorporates graph structure linearly, significantly reducing the pattern count and computational load. This efficiency, coupled with the ability to integrate topological information, makes \(\mvDEG\) a superior choice for complex data analysis, especially in multichannel scenarios, as in Sec.~\ref{sec:time}.

\subsection{Correlated Noise Analysis}\label{sub:corr}

In this section, we explore the behaviour of multivariate methods on correlated noise. The experiment involves generating four sets of time series, each with distinct correlation structures: (1) Four Uncorrelated WGN Series, (2) Two Correlated and Two Uncorrelated Series (a mix of correlated and uncorrelated series), (3) Two Pairs of Correlated Series (Two distinct pairs, each internally correlated), (4) Three Correlated and One Uncorrelated Series and (5) Four Correlated Series.

For each set of time series, we computed both \(\mvDE\) and \(\mvDEG\) using \(N=500\), \(m=4\), and \(c=6\). The mean entropy values and their SD were determined over $40$ simulations. In these computations, the underlying graph \(\G\) was based on the theoretical correlation matrix. When using the estimated correlation matrix \(\hat{\G}\), derived from the data, the entropy values obtained for each simulation were remarkably consistent, indicating the robustness of \(\mvDEG\) in capturing the dynamics of the time series even with small variations in the graph structure.
 
Additionally, our focus on shorter time series (\(N=500\)) demonstrates $\mvDEG$'s robustness with limited data, a frequent challenge where traditional methods struggle. This underscores the method's adaptability to scenarios with constrained data, as exemplified in Sec.~\ref{sec:real}. Comparable results are observed with larger sample sizes.

\begin{figure}
	\centering
	\includegraphics[width=.45\textwidth]{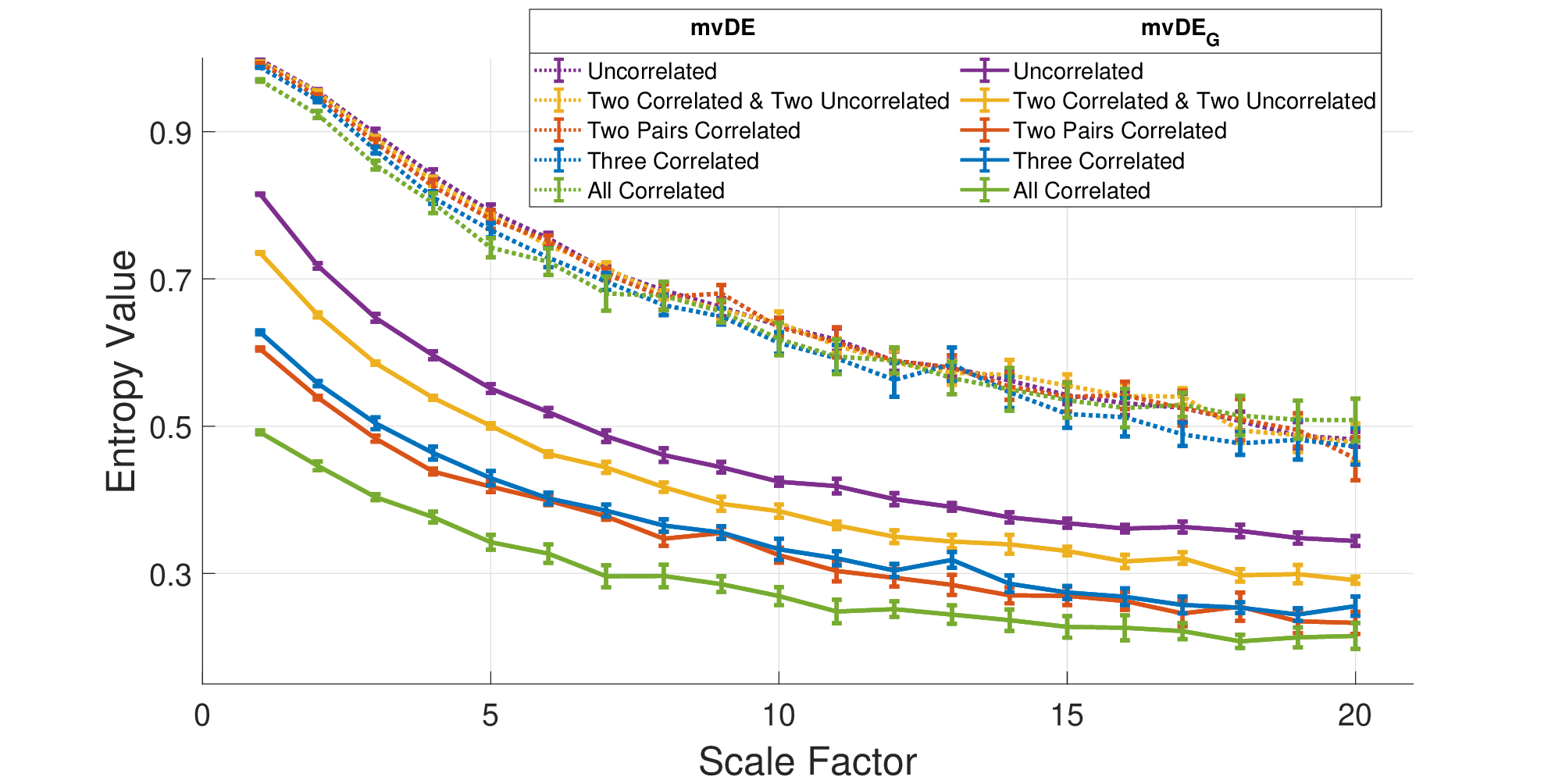}
	\caption{Correlated noise entropy across various structures, ranging from uncorrelated to fully connected across scale factors.}
	\label{fig:Correlated}
\end{figure}

The results, as shown in Fig.~\ref{fig:Correlated}, reveal interesting insights into the dynamics of correlated noise. In both \(\mvDE\) and \(\mvDEG\) methods, the entropy values for WGN decrease as the scale factor increases, indicating a reduction in complexity. However, a key observation in mvDE is the overlapping of entropy values across all sets of time series and scales, suggesting that mvDE struggles to differentiate between various degrees of correlation and uncorrelated dynamics.

In contrast, our proposed method, $\mvDEG$, not only presents distinct mean entropy values for each set of time series but also avoids the overlapping of standard deviations. This distinct separation in the entropy values and their variances demonstrates that $\mvDEG$ is more adept at characterizing the dynamics of time series with varying degrees of correlation. This capability makes $\mvDEG$ a more robust tool for analysing complex multivariate time series data, especially in scenarios where understanding the interplay between correlation and signal complexity is crucial.

\begin{figure}
	\centering
	\includegraphics[width=.45\textwidth]{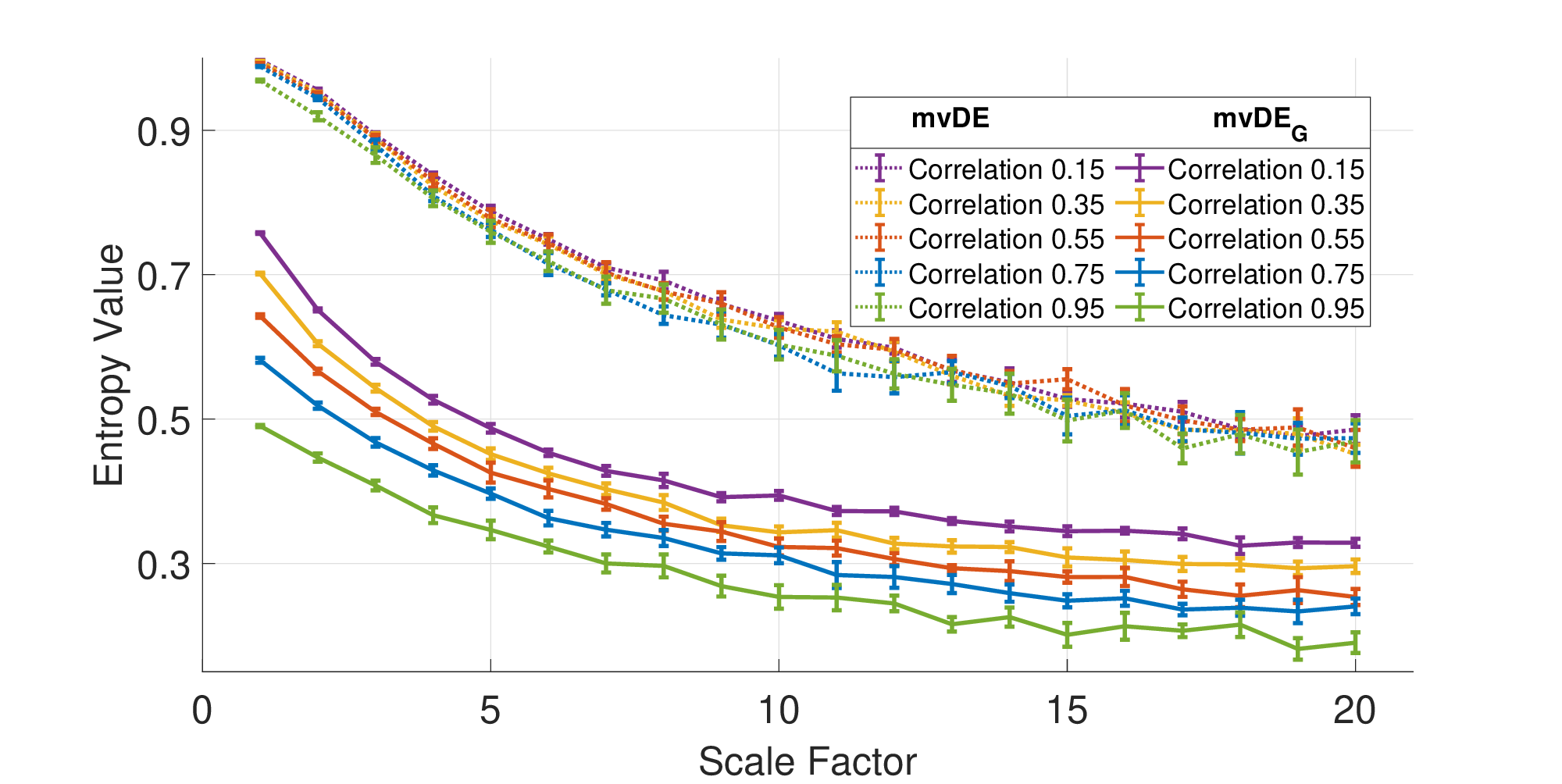}
	\caption{Entropy for varying degrees of correlated noise.}
	\label{fig:degree}
\end{figure}

Further, we explore the dynamics of correlated noise with varying degrees of correlation. For this purpose, we generate three time series for each set, with correlation values set at 0.95, 0.75, 0.55, 0.35, and 0.15. The parameters used are \(m=4\), \(c=6\), and the analysis is based on 40 realizations. The results, depicted in Fig.~\ref{fig:degree}, include both the mean and SD of the entropy values.

The analysis reveals that while both methods show a decrease in entropy values with increasing scale factor, the classical $\mvDE$ fails to distinguish the dynamics influenced by varying degrees of correlation, as evidenced by the overlapping SDs. In contrast, $\mvDEG$ successfully differentiates between these dynamics, with minimal overlap in the SDs across different correlation levels. This distinction is not only pivotal for understanding complex signal interactions but also highlights the superior computational efficiency of $\mvDEG$, especially for shorter time series. The following section will further elaborate on the computational advantages of our proposed method.

\section{Computational Time Comparison}\label{sec:time}
This section compares the computational efficiency of classical \(\mvDE\) and \(\mvDEG\). To assess the computational performance, we utilized uncorrelated multivariate WGN time series of varying lengths, ranging from $100$ to $10,000$ sample points, across $10$ channels. The experiments were conducted on a PC running MATLAB R2023b, equipped with 32.0 GB RAM and an Intel(R) Core(TM) i7-10610U CPU @ 1.80GHz, 2.30 GHz processor.

\begin{figure}
	\centering
	\includegraphics[width=.45\textwidth]{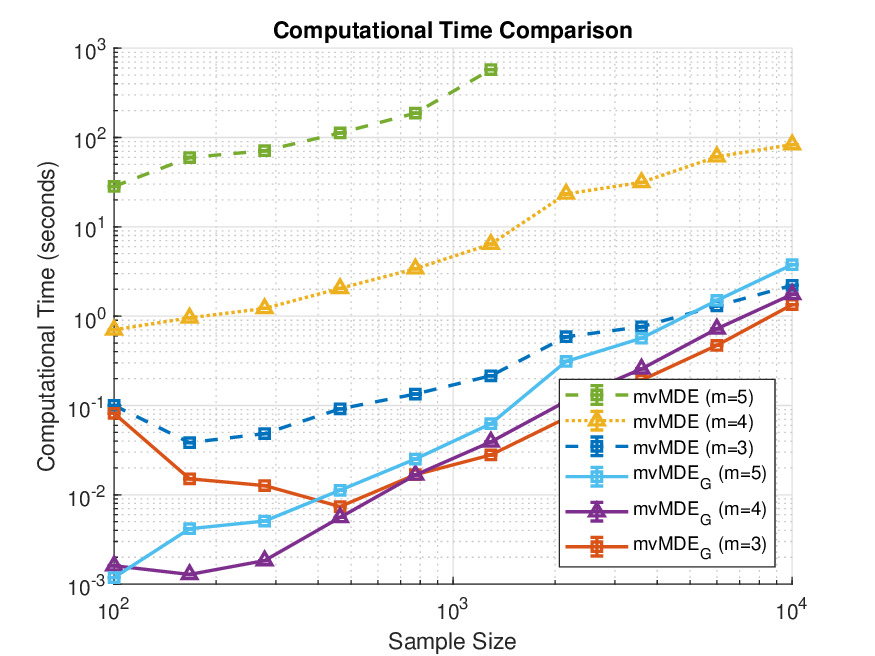} % Replace with your figure file
	\caption{Computational times comparison for \(\mvDEG\) against classical \(\mvDE\), for different data lengths.}
	\label{fig:time_comparison}
\end{figure}

The classical \(\mvDE\) method exhibited significant limitations in processing large channels. For a dataset with \(\sim2,000\) sample points, \(m=5\) and \(p=8\) channels, MATLAB encountered an overflow RAM error. This limitation is problematic in scenarios requiring the analysis of high-dimensional data or multiple channels. In contrast, \(\mvDEG\) showcased computational efficiency under analogous conditions. It successfully processed data without encountering the memory constraints that impeded the classical method. As depicted in Fig.~\ref{fig:time_comparison}, \(\mvDEG\) consistently outperformed \(\mvDE\) in terms of computational time, especially in datasets with a larger number of sample points and channels.

The classical \(\mvDE\) method computes $(N-m+1)\binom{mp}{m}$ dispersion patterns. This combinatorial approach, while making \(\mvDE\) robust to very short data sequences by generating a multitude of patterns even from a limited time series, also introduces a computational complexity. The quantity of such dispersion patterns highlights the exponential growth in the number with increasing $m$ and $p$, leading to increases in both computational time and resource usage.

The superior performance of \(\mvDEG\) is primarily attributed to its algorithmic design. The computational requirement of \(\mvDEG\) is significantly lower, as it is limited to processing at most $(N-m)p$ dispersion patterns and also the algorithm optimizes the computation of $Np \times Np$ matrix powers. This difference in computational requirements underscores the advanced algorithmic design of \(\mvDEG\), making it a preferable choice for applications involving large-scale multivariate time series analysis even if the results are similar. For instance, in a dataset with $N=2000$, $p=8$, and $m=5$, the classical \(\mvDE\) method needs to compute approximately $1.3 \times 10^9$ dispersion patterns. In contrast, \(\mvDEG\) requires only about $1.6 \times 10^3$ patterns, demonstrating a substantial reduction in computational complexity. 

This efficiency, combined with its robustness for short time series, positions \(\mvDEG\) as a superior tool in time series analysis.

\section{Real-World Applications}\label{sec:real}

This section demonstrates the application of $\mvDEG$ to meteorological data and industrial two-phase flow systems. The complex nature of weather data~\cite{fraedrich1986estimating} and the dynamic characteristics of industrial flows~\cite{tan2023combinational} provide ideal test to showcase the efficacy and adaptability of $\mvDEG$.

\subsection{Weather Data on Ground Station Graphs}
\subsubsection{Data Description}
We analysed temperature data from Brittany ground stations (January 2014), provided by the French national meteorological service~\cite{girault2015stationary}. A graph was constructed where vertices represent stations, and edges' weights are based on the Gaussian kernel of the Euclidean distances between stations. This setup allows for a detailed examination of spatial interactions and temperature variability across the region. The weighting formula between two vertices \(i\) and \(j\) is as follows:
\begin{equation}
	\A_{i j}= 
	\begin{cases}
		\exp \left(\frac{-d(i, j)^2}{2 \sigma_1^2}\right) & \text{if } d(i, j) \leq \sigma_2 \\
		0 & \text{otherwise.}
	\end{cases}
	\label{eq:weight}
\end{equation}
where \(d(i, j)\) denotes the Euclidean distance between stations \(i\) and \(j\), with parameters \(\sigma_1^2=5.1^8\) and \(\sigma_2=10^5\), as per~\cite{girault2015stationary}.

\begin{figure}
	\centering
	\includegraphics[width=.45\textwidth]{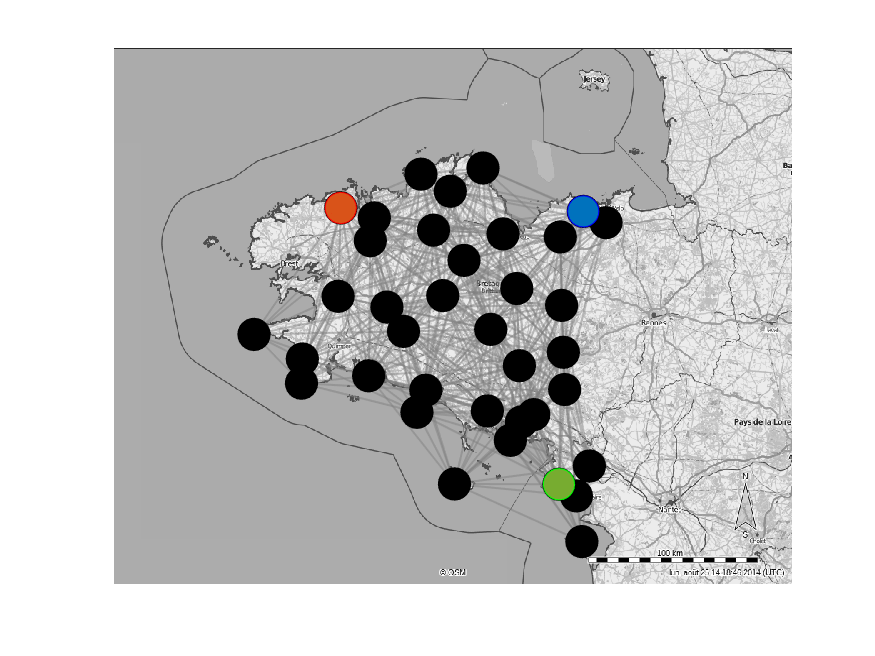}
	\caption{Geographical distribution of Brittany ground station, with selected cities highlighted for detailed analysis.}
	\label{fig:map}
\end{figure}

\begin{figure}
	\centering
	\includegraphics[width=.85\columnwidth]{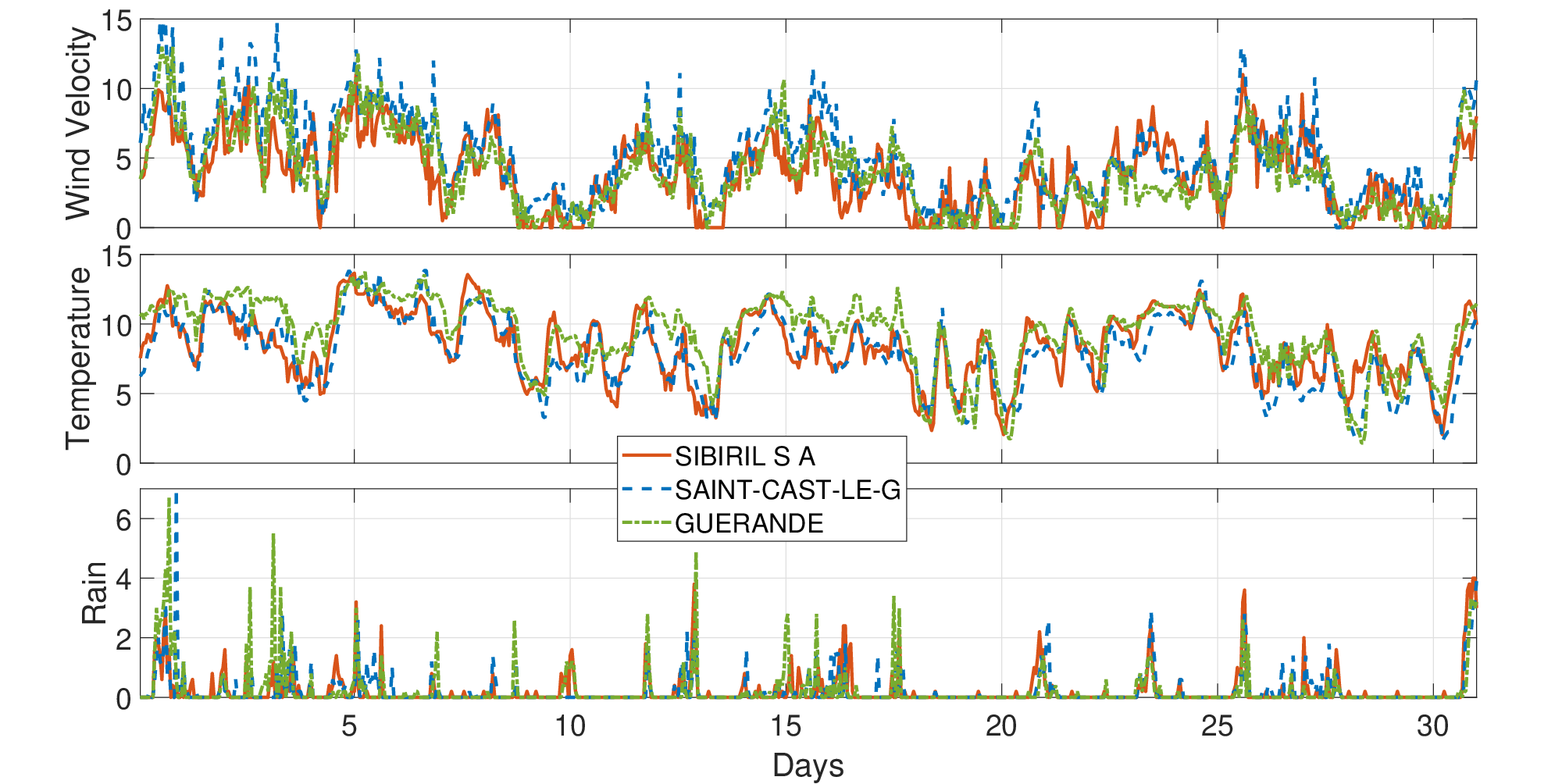}
	\caption{Representative time series for wind velocity, temperature, and rainfall recorded at the Brittany stations over a 30-day period, as identified in Fig.~\ref{fig:map}.}
	\label{fig:timemap}
\end{figure}

\begin{figure}
	\centering
	\includegraphics[width=.4\textwidth]{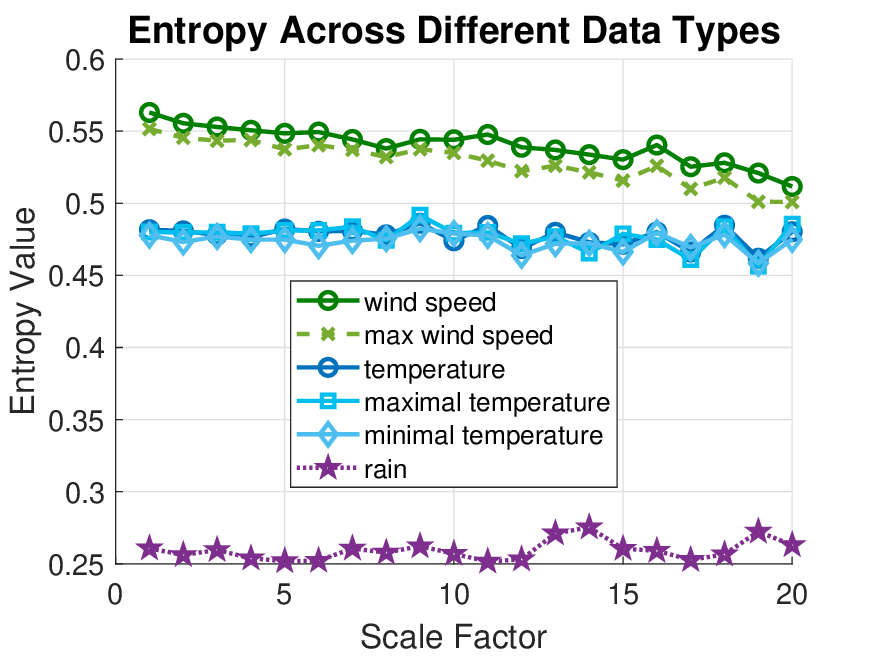}
	\caption{Entropy of meteorological parameters: wind speed, maximum wind speed, temperature, and rainfall.}
	\label{fig:temperature}
\end{figure}

Fig.~\ref{fig:map} depicts Brittany's map with sensor locations as black nodes. The edges are computed using Eq.~\ref{eq:weight}. Three cities (vertices) are marked in different colours (red, green, and blue), and their corresponding time series for wind velocity, temperature, and rainfall are shown in Fig.~\ref{fig:timemap}. Each time series consists of \(744\) samples, representing hourly data for \(31\) days.
 The entropy values for \(\mvDEG\) were computed using the distance graph for \(m=4\), \(c=6\), and the results are presented in Fig.~\ref{fig:temperature}. 
\subsubsection{Results and Discussion}
The entropy analysis distinctly identifies patterns in meteorological data. Rainfall data (purple) consistently show low entropy values across scales, akin to the characteristics of \(1/f\) noise, indicating a consistent, but low, level of complexity. In comparison, temperature data (blue) display slightly higher, yet stable entropy values, pointing to a different pattern of complexity. Wind data (mean and max wind speed) demonstrate the greatest complexity, with a decrease at lower scales.

\(\mvDEG\) effectively captures the varying characteristics of meteorological elements, especially the complex and fluctuating patterns of wind~\cite{laib2018multifractal}. Influenced by factors like local terrain and atmospheric conditions, wind patterns pose significant challenges for accurate prediction and modelling. In contrast, temperature changes tend to be more predictable, following daily and seasonal cycles. Rainfall, subject to variability, generally follows recognizable weather patterns, showing less sudden changes compared to the more unpredictable nature of wind.

The entropy analysis distinctly differentiates meteorological elements, with higher entropy in wind data emphasizing the challenges in its prediction and modelling – a critical factor for industries reliant on accurate wind forecasts. In contrast, the more stable and lower entropy patterns in rainfall data underscore $\mvDEG$'s capacity to effectively delineate the varied complexities of meteorological phenomena.

\subsubsection{Comparison with other state-of-the-art methods}

Despite the relatively small size of our dataset, consisting of data from \(37\) cities each with \(744\) samples, the computational demand of the classical \(\mvDE\) is substantial. Specifically, it requires the processing of approximately \( 3.4\times 10^{11} \) patterns. This high computational load, largely due to each city being treated as a separate channel, renders the classical \(\mvDE\) impractical for execution on standard personal computers because it leads to memory errors. . In contrast, our \(\mvDEG\) method drastically reduces this requirement to only about \(2.7 \times 10^4\) patterns. This significant decrease in computational complexity not only makes \(\mvDEG\) feasible for processing large datasets on typical computing setups but also extends its applicability to a broader range of real-world scenarios.

\begin{figure}[h]
	\centering
	\begin{subfigure}{0.24\textwidth}
		\centering
		\includegraphics[width=\textwidth]{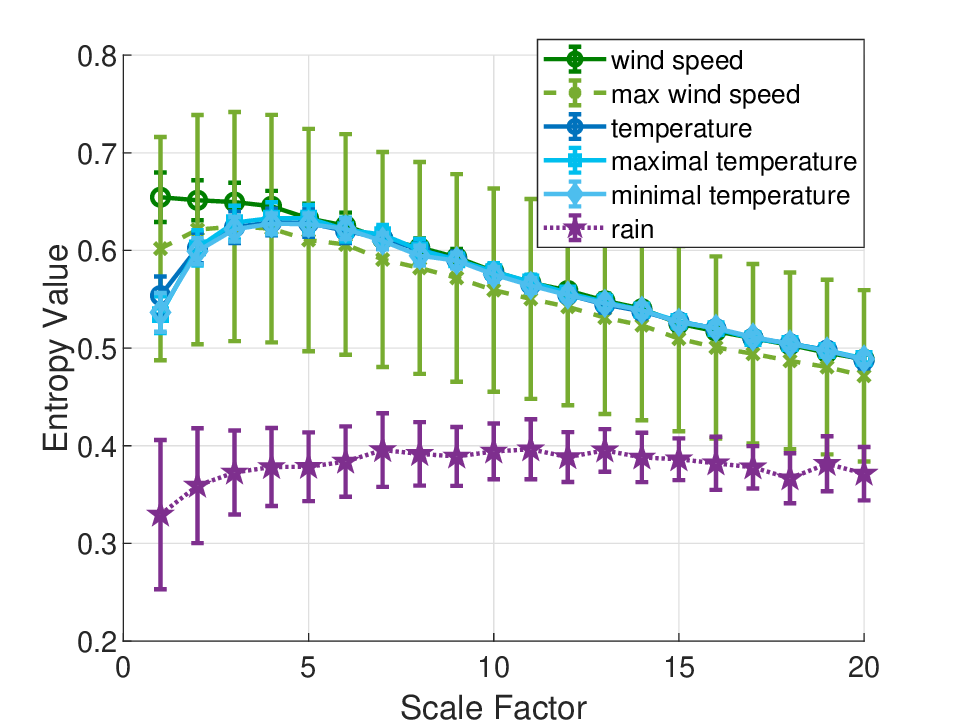}
		\caption{}
		\label{fig:mvde_empty1}
	\end{subfigure}
	\begin{subfigure}{0.24\textwidth}
		\centering
		\includegraphics[width=\textwidth]{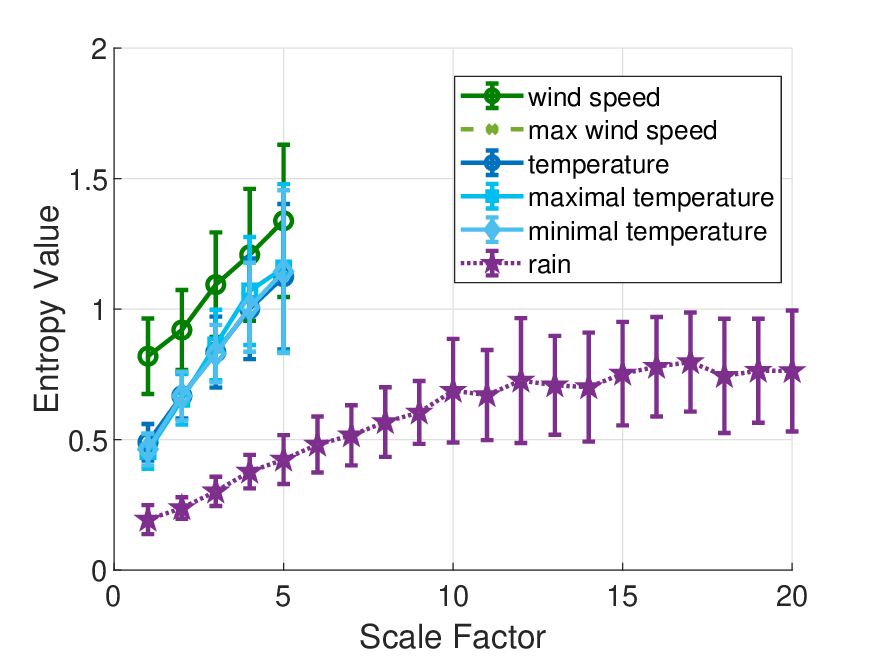}
		\caption{}
		\label{fig:mvdeg_corr2}
	\end{subfigure}
	\caption{Comparison of entropy analysis of weather data with state-of-the-art methods: (a) Dispersion Entropy~\cite{Rostaghi2016}, (b) Sample Entropy~\cite{azami2018coarse}.}
	\label{fig:temperature2}
\end{figure}

Due to the impracticality of applying classical \(\mvDE\) to our dataset, we resorted to using the univariate version of the Multiscale Dispersion Entropy and Sample Entropy for analysing the weather data. We use Sample Entropy due to its widespread application in data analysis~\cite{huo2020entropy}.
This approach limits our analysis to the temporal information of each sensor, disregarding the spatial relationships between them. The mean entropy results from both methods are presented in Fig.~\ref{fig:temperature2}. While these methods are effective in detecting the dynamics and yield lower entropy values for rainfall data, they fail to distinguish the complex dynamics between temperature and wind data. Additionally, Sample Entropy becomes undefined~\cite{azami2018coarse} for scale factors~\(>6\).

\subsection{Two-phase Flow Data}

\subsubsection{Data Description}
We analyse data from two-phase flow experiments conducted with Electrical Resistance Tomography (ERT), a technique that leverages the varying conductivities of mediums to reconstruct the conductivity distribution within a pipe. This is achieved by measuring boundary voltages between electrodes, with a constant electrical current of 50 kHz serving as the excitation signal and a data acquisition rate of 120 frames per second~\cite{Tan2013, tan2023combinational}. The $16-$electrode ERT setup captures $16\times 13=208$ voltage data points. To manage the redundancy in these measurements, feature vectors $V_{Ri}$ are extracted from each electrode, reducing the dimensionality of the data. This process involves averaging the voltage measurements relative to a baseline condition (when the pipe is full of water), calculated as $V_{Ri}=\frac{1}{13}\sum_{j=1}^{13}(V_{ij}-V_{i{j_0}})/V_{i{j_0}}$, where $V_{ij}$ is the measured voltage value and $V_{i{j_0}}$ is the baseline voltage for electrode $i$~\cite{Tan2013}. In the experiments, gas and water were mixed at velocities ranging from \(0.4 \, \text{m/s}\) to \(2.9 \, \text{m/s}\) for water and \(0.06 \, \text{m/s}\) to \(5.64 \, \text{m/s}\) for gas. This produced diverse flow patterns, which includes \(105\) distinct experiments conducted under various gas and water flow rates, each generated around \(1400\) time samples.

Two-phase flow experiments usually require data dimensionality reduction for computational analysis, often by averaging feature vectors \( V_{Ri} \) from a 16-electrode ERT system into four aggregated time series, as per \cite{Tan2013}. In contrast, \(\mvDEG\) processes the entire dataset directly, eliminating the need for initial data reduction. This capability allows \(\mvDEG\) to efficiently handle the full dataset, including raw, unfiltered data, demonstrating its robustness against noise.

\subsubsection{Results and Discussion}
\begin{figure}
	\centering
	\includegraphics[width=.45\textwidth]{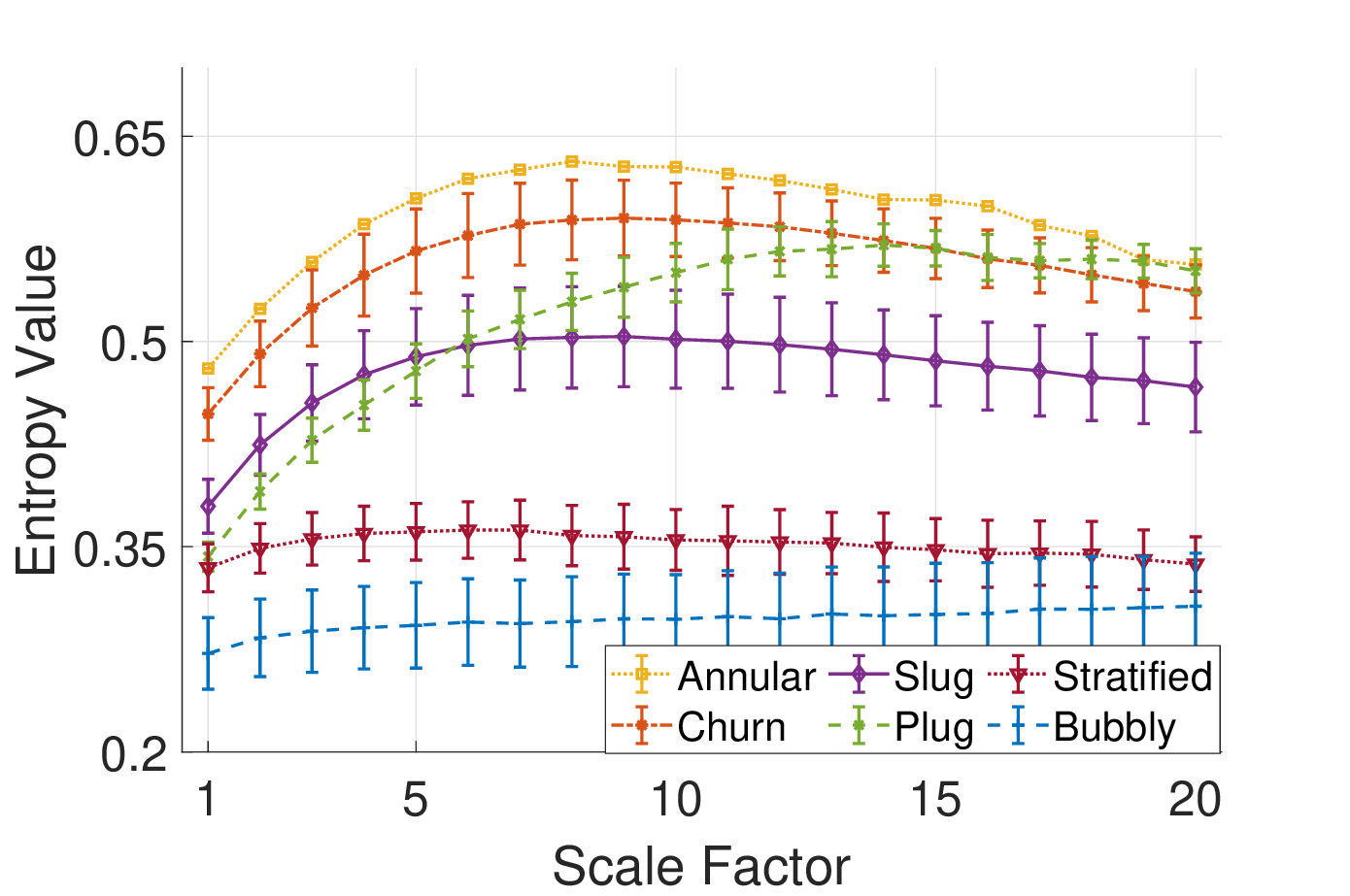} % Replace with your figure file
	\caption{Entropy values for the two-phase flow data.}
	\label{fig:phase-flow}
\end{figure}

For our analysis, we applied the \(\mvDEG\) technique, selecting \(m = 4\) and \(c = 6\). This choice of the graph was based on the fact that each \(V_{Ri}\) represents an unweighted average of nodes, suggesting the use of a complete graph as the underlying structure. We found consistent results across various underlying graphs derived from the data, as well as with different \(m\) and \(c\) values, and even when analysing cleaned data. The entropy values obtained for the six flow patterns were plotted against the scale factor in Fig.~\ref{fig:phase-flow}.

A key element of our \(\mvDEG\) analysis is assessing entropy variations across scales, offering vital insights into the dynamics of each flow pattern. For example, the consistent entropy profile in the bubbly flow regime suggests uniform complexity and predictability, while the entropy shifts in slug and annular flows reflect changes in system complexity and dynamics at various scales. These fluctuations are linked to alterations in flow structure and phase interactions with changing scale factors.

\subsubsection*{Bubbly Flow}
The bubbly flow regime is marked by the lowest entropy values, indicative of a structured and predictable system. This pattern results from small bubbles' limited impact on selected electrodes, intensified by gravity's tendency to concentrate bubbles at the top of the electrode array. Therefore, the sensors exhibit regular readings, leading to lower entropy values~\cite{Tan2013}. As the scale factor increases, the influence of these bubbles becomes more noticeable, causing slight entropy variations. Our \(\mvDEG\) method effectively captures this nuanced impact, especially in electrodes more directly exposed to the gas phase. The overall flattened entropy profile across scales in this regime reflects its stable and predictable nature.
%
%\subsubsection*{Bubbly Flow} In the bubbly flow regime, we observe the lowest entropy values, indicative of a more structured and predictable system. This pattern arises primarily due to the limited impact of small bubbles on a select number of electrodes, an effect that is intensified by gravity, which tends to concentrate bubbles towards the top of the electrodes array. As a result, the sensor, represented by the vertices (electrodes) in our analysis, exhibit regular or nearly constant readings, leading to lower entropy values~\cite{Tan2013}. Interestingly, as the scale factor increases, the influence of the bubbles becomes more pronounced, causing slight variations and a marginal increase in entropy values at larger scales. This pattern is particularly evident in electrodes that are more directly exposed to the gas phase, where fluctuations are relatively larger. Our \(\mvDEG\) method adeptly captures these nuances, effectively highlighting the subtle yet significant impact of bubbles on the surrounding fluid medium.
%
%The flattened entropy profile over scales in the bubbly flow regime reflects a consistent dynamic behaviour, underscoring its structured and predictable nature.

\subsubsection*{Stratified Flow} Stratified flow is characterized by a low and stable entropy profile, indicative of a system with reduced complexity. This flow regime, marked by a clear separation between the gas and liquid phases, exhibits minimal fluctuations across all \(16\) electrodes. The \(\mvDEG\) method's ability to delineate between closely related flow patterns, particularly stratified and bubbly flows, is evident in lower scales where traditional methods might falter~\cite{Tan2013}. This capability to separate stratified flow from other flow patterns, even at lower scales, highlights the robustness and precision of the \(\mvDEG\) method in analysing complex two-phase flow dynamics.

\subsubsection*{Slug Flow}
Slug flow displays distinctive entropy values, with higher values at lower scales, indicative of the pronounced influence of large gas bubbles. This complexity is most evident in the initial scales, as the large bubbles' dynamics generate unique patterns. As shown in Fig.~\ref{fig:phase-flow}, a decreasing trend in entropy values emerges beyond scale factor 10, reflecting a change in flow complexity across scales. \(\mvDEG\) adeptly differentiates slug flow from other flow regimes, particularly at lower scales, underscoring its capacity to capture the intricate dynamics unique to this type of flow. The entropy profile, with its initial rise and subsequent fall, mirrors the varying influence of large gas bubbles within slug flow, revealing their significant impact on the flow's overall complexity.

\subsubsection*{Plug Flow} In contrast to slug flow, plug flow exhibits lower entropy values, particularly at smaller scales. This is attributed to the complex dynamics of smaller bubbles, which produce more variation and an increase in entropy values as the scale factor increases. However, a notable transition occurs around scale factor 10, where the entropy values for plug flow diverge significantly from those of slug flow. This divergence is indicative of the transition from plug to slug flow, characterized by the coalescence of smaller bubbles into larger ones. The transition from plug flow to slug flow is distinct from bubbly flow across all scales considered. With increasing gas velocity, a complexity transition is observed at larger scales, and the transition scale tends to shift towards lower scales. This shift, as depicted in Fig.~\ref{fig:phase-flow}, suggests that scale 10 is particularly sensitive to the transition from plug flow to slug flow, as also supported by similar findings in the literature (see~\cite{Tan2013}).

\subsubsection*{Churn Flow} Churn flow, occupying an intermediate position between slug and annular flows, is indicative of the system's increasing complexity. This flow regime is characterized by entropy values that surpass those of slug flow but remain below the peaks of annular flow, embodying a transitional state with elements of both adjacent patterns. \cite{taitel1980modelling} considered churn flow as a developing slug flow, while \cite{hewitt2013annular} described it as intermediate between slug and annular flows. This characterization aligns with our entropy analysis, where churn flow's entropy values are consistently higher than slug flow but lower than annular flow across all scale factors. The churn regime, known for its chaotic nature and complexity~\cite{ansari1994comprehensive}, is accurately reflected in the high entropy values obtained from our analysis. This phase serves as a critical indicator of the system's transition towards the more chaotic nature characteristic of annular flow.

\subsubsection*{Annular Flow}
Annular flow is characterized by the highest entropy values, indicative of its complex and turbulent nature. This regime involves a rapidly moving gas core surrounded by a liquid film with embedded gas bubbles and droplets, contributing to its high complexity~\cite{tan2023combinational}. In horizontal annular flow, gravity causes an uneven liquid film distribution, creating surface waves. These waves, disrupted by gas shear forces, produce droplets in the gas core, adding to the regime's intricacy~\cite{Tan2013}. \(\mvDEG\) effectively captures the dynamic interplay between these elements, highlighting the chaotic characteristics of annular flow. The entropy values initially rise and then stabilize at higher scales, reflecting a transition from variable to more consistent complexity within this regime.

Our \(\mvDEG\) analysis captures the unique entropy profiles for each flow pattern in two-phase flow, providing insights into the underlying dynamics. The method's sensitivity to nuanced dynamics is evident in the clear separation of entropy values for each regime, particularly at lower scales. Additionally, \(\mvDEG\) allows us to utilize the total information from all 16 electrodes without the need to reduce them to a smaller set of vertices, thus preserving the richness of the data. This capability contrasts with the \(\mvDE\) approach, where attempting to use all \(16\) electrodes would lead to significant computational difficulties.  Indeed, even with only 8 channels and a similar size of the sample of two-phase flow data, the state-of-the-art methods, as detailed in Sec.~\ref{sec:time}, have been observed to produce out-of-memory errors due to their computational demands. In this way, \(\mvDEG\) not only offers more detailed insights but also ensures computational feasibility when handling extensive sensor data.

\section{Conclusion}
\label{sec:conclusion}

This paper presented the Multivariate Multiscale Graph-based Dispersion Entropy (\(\mvDEG\)), a new method for analysing multivariate time series on graphs and networks. \(\mvDEG\) uniquely combines temporal dynamics with topological structure, enhancing analysis beyond traditional entropy methods.

\(\mvDEG\) has proven effective in differentiating complex patterns in multivariate time series, excelling particularly with short series where conventional methods falter. Its application to real-world datasets, including weather and two-phase flow data, highlights its versatility and robustness in practical scenarios.

A key innovation of \(\mvDEG\) is its computational efficiency, a significant leap forward from existing methods. By utilizing matrix properties and the Kronecker product, \(\mvDEG\) transitions computational time growth from the exponential rates typical of classical methods to a linear growth with respect to the number of vertices or nodes, making it apt for extensive and real-time analyses.

In summary, \(\mvDEG\) stands as a powerful tool in multivariate time series analysis, adept at unraveling complex spatial-temporal interplays in various fields.

\bibliographystyle{abbrv}
\bibliography{references_reduced} % This is your .bib file's name without the extension

\begin{thebibliography}{10}

\bibitem{Ahmed2011}
M.~U. Ahmed and D.~P. Mandic.
\newblock Multivariate multiscale entropy.
\newblock {\em Phys. Rev. E - Stat. Nonlinear, Soft Matter Phys.},
  84(6):061918, Dec 2011.

\bibitem{ahmed2011multivariate}
M.~U. Ahmed and D.~P. Mandic.
\newblock Multivariate multiscale entropy analysis.
\newblock {\em IEEE Signal Processing Letters}, 19(2):91--94, 2011.

\bibitem{ansari1994comprehensive}
A.~Ansari, N.~Sylvester, C.~Sarica, O.~Shoham, and J.~Brill.
\newblock A comprehensive mechanistic model for upward two-phase flow in
  wellbores.
\newblock {\em SPE Production \& Facilities}, 9(02):143--151, 1994.

\bibitem{azami2018coarse}
H.~Azami and J.~Escudero.
\newblock Coarse-graining approaches in univariate multiscale sample and
  dispersion entropy.
\newblock {\em Entropy}, 20(2):138, 2018.

\bibitem{Azami2019}
H.~Azami, A.~Fernández, and J.~Escudero.
\newblock Multivariate multiscale dispersion entropy of biomedical time series.
\newblock {\em Entropy}, 21(9):1--21, 2019.

\bibitem{boccaletti2006complex}
S.~Boccaletti and et~al.
\newblock Complex networks: Structure and dynamics.
\newblock {\em Physics Reports}, 424(4-5):175--308, 2006.

\bibitem{bondy1982graph}
J.~A. Bondy.
\newblock {\em Theory with Applications}.
\newblock 1982.

\bibitem{chai2019enhanced}
Z.~Chai and C.~Zhao.
\newblock Enhanced random forest with concurrent analysis of static and dynamic
  nodes for industrial fault classification.
\newblock {\em IEEE Trans Industr Inform}, 16(1):54--66, 2019.

\bibitem{costa2002multiscale}
M.~Costa, A.~L. Goldberger, and C.-K. Peng.
\newblock Multiscale entropy analysis of complex physiologic time series.
\newblock {\em Physical Review Letters}, 89(6):068102, 2002.

\bibitem{dang2019novel}
W.~Dang and et~al.
\newblock A novel deep learning framework for industrial multiphase flow
  characterization.
\newblock {\em IEEE Trans Industr Inform}, 15(11):5954--5962, 2019.

\bibitem{dong2024information}
K.~Dong and D.~Li.
\newblock An information theoretic approach to analyze irregularity of graph
  signals and network topological changes based on bubble entropy.
\newblock {\em Fluctuation and Noise Letters}, 2024.

\bibitem{fabila2022multivariate}
J.~S. Fabila-Carrasco, C.~Tan, and J.~Escudero.
\newblock Multivariate permutation entropy, a cartesian graph product approach.
\newblock In {\em 2022 30th European Signal Processing Conference (EUSIPCO)},
  pages 2081--2085. IEEE, 2022.

\bibitem{Fabila22}
J.~S. Fabila-Carrasco, C.~Tan, and J.~Escudero.
\newblock Permutation entropy for graph signals.
\newblock {\em IEEE Trans. Signal Inf. Process. over Networks}, 8:288--300,
  2022.

\bibitem{Fabila23}
J.~S. Fabila-Carrasco, C.~Tan, and J.~Escudero.
\newblock Dispersion entropy for graph signals.
\newblock {\em Chaos, Solitons \& Fractals}, 175:113977, Oct 2023.

\bibitem{fan2023graph}
Q.~Fan and et~al.
\newblock Graph multi-scale permutation entropy for bearing fault diagnosis.
\newblock {\em Sensors}, 24(1):56, 2023.

\bibitem{fraedrich1986estimating}
K.~Fraedrich.
\newblock Estimating the dimensions of weather and climate attractors.
\newblock {\em Journal of Atmospheric Sciences}, 43(5):419--432, 1986.

\bibitem{girault2015stationary}
B.~Girault.
\newblock Stationary graph signals using an isometric graph translation.
\newblock In {\em 2015 EUSIPCO}, pages 1516--1520. IEEE, 2015.

\bibitem{graham2018kronecker}
A.~Graham.
\newblock {\em Kronecker Products and Matrix Calculus with Applications}.
\newblock Courier Dover Publications, 2018.

\bibitem{hewitt2013annular}
G.~Hewitt.
\newblock {\em Annular two-phase flow}.
\newblock Elsevier, 2013.

\bibitem{humeau2016multivariate}
A.~Humeau-Heurtier.
\newblock Multivariate generalized multiscale entropy analysis.
\newblock {\em Entropy}, 18(11):411, 2016.

\bibitem{huo2020entropy}
Z.~Huo, M.~Mart{\'\i}nez-Garc{\'\i}a, Y.~Zhang, R.~Yan, and L.~Shu.
\newblock Entropy measures in machine fault diagnosis: Insights and
  applications.
\newblock {\em IEEE Transactions on Instrumentation and Measurement},
  69(6):2607--2620, 2020.

\bibitem{jhaveri2021fault}
R.~H. Jhaveri, S.~V. Ramani, G.~Srivastava, T.~R. Gadekallu, and V.~Aggarwal.
\newblock Fault-resilience for bandwidth management in industrial
  software-defined networks.
\newblock {\em IEEE Transactions on Network Science and Engineering},
  8(4):3129--3139, 2021.

\bibitem{laib2018multifractal}
M.~Laib, J.~Golay, L.~Telesca, and M.~Kanevski.
\newblock Multifractal analysis of the time series of daily means of wind speed
  in complex regions.
\newblock {\em Chaos, Solitons \& Fractals}, 109:118--127, 2018.

\bibitem{ortega2018graph}
A.~Ortega, P.~Frossard, J.~Kova{\v{c}}evi{\'c}, J.~M. Moura, and
  P.~Vandergheynst.
\newblock Graph signal processing: Overview, challenges, and applications.
\newblock {\em Proceedings of the IEEE}, 106(5):808--828, 2018.

\bibitem{prasse2020network}
B.~Prasse and P.~Van~Mieghem.
\newblock Network reconstruction and prediction of epidemic outbreaks for
  general group-based compartmental epidemic models.
\newblock {\em IEEE Transactions on Network Science and Engineering},
  7(4):2755--2764, 2020.

\bibitem{qiu2017time}
K.~Qiu, X.~Mao, X.~Shen, X.~Wang, T.~Li, and Y.~Gu.
\newblock Time-varying graph signal reconstruction.
\newblock {\em IEEE Journal of Selected Topics in Signal Processing},
  11(6):870--883, 2017.

\bibitem{Rostaghi2016}
M.~Rostaghi and H.~Azami.
\newblock Dispersion entropy: A measure for time-series analysis.
\newblock {\em IEEE Signal Process. Lett.}, 23(5):610--614, May 2016.

\bibitem{taitel1980modelling}
Y.~Taitel, D.~Barnea, and A.~Dukler.
\newblock Modelling flow pattern transitions for steady upward gas-liquid flow
  in vertical tubes.
\newblock {\em AIChE Journal}, 26(3):345--354, 1980.

\bibitem{tan2023combinational}
C.~Tan, H.~Jia, G.~Liang, X.~Wang, W.~Niu, and F.~Dong.
\newblock Combinational multi-modality tomography system for industrial
  multiphase flow imaging.
\newblock {\em IEEE Trans Instrum Meas}, 2023.

\bibitem{Tan2013}
C.~Tan, J.~Zhao, and F.~Dong.
\newblock Gas-water two-phase flow characterization with electrical resistance
  tomography and multivariate entropy analysis.
\newblock {\em ISA Transactions}, 55, 2015.

\bibitem{wang2021variational}
X.~Wang and et~al.
\newblock Variational embedding multiscale diversity entropy for fault
  diagnosis of large-scale machinery.
\newblock {\em IEEE Transactions on Industrial Electronics}, 69(3):3109--3119,
  2021.

\bibitem{wei2019optimal}
Z.~Wei, A.~Pagani, G.~Fu, I.~Guymer, W.~Chen, J.~McCann, and W.~Guo.
\newblock Optimal sampling of water distribution network dynamics using graph
  fourier transform.
\newblock {\em IEEE Transactions on Network Science and Engineering},
  7(3):1570--1582, 2019.

\bibitem{xia2021graph}
F.~Xia and et~al.
\newblock Graph learning: A survey.
\newblock {\em IEEE Transactions on Artificial Intelligence}, 2(2):109--127,
  2021.

\bibitem{zhou2020edm}
R.~Zhou and et~al.
\newblock Edm-fuzzy: An euclidean distance-based multiscale fuzzy entropy
  technology for diagnosing faults of industrial systems.
\newblock {\em IEEE Trans Industr Inform}, 17(6), 2020.

\end{thebibliography}

\end{document}